# BI-OBJECTIVE OPTIMIZATION OF A VRP PROBLEM APPLIED TO URBAN SOLID WASTE COLLECTION THROUGH A MODEL THAT INCLUDES THE VISUAL ATTRACTION OF ROUTES

# OPTIMIZACIÓN BIOBJETIVO DE UN PROBLEMA VRP APLICADO A LA RECOLECCIÓN DE RESIDUOS SÓLIDOS URBANOS MEDIANTE UN MODELO QUE INCLUYE LA ATRACCIÓN VISUAL DE LAS RUTAS

# OTIMIZAÇÃO BI-OBJETIVO DE UM PROBLEMA VRP APLICADA À COLETA DE RESÍDUOS SÓLIDOS URBANOS POR MEIO DE UM MODELO QUE INCLUI A ATRAÇÃO VISUAL DAS ROTAS

Diego Gabriel Rossit, INMABB Departamento de Ingeniería - Universidad Nacional del Sur CONICET, diego.rossit@uns.edu.ar

Adrián Andrés Toncovich, Departamento de Ingeniería - Universidad Nacional del Sur, atoncovi@uns.edu.ar

**Abstract**

The compactness of routes in distribution plans is a criterion that has not been sufficiently explored in the literature related to logistics distribution but has shown to have a significant impact on the practical implementation of routing plans, for example in solid waste collection problems. In this regard, this article presents a bi-objective model to optimize the vehicle routing problem with time constraints, considering the minimization of travel times and the compactness of routes. Experimental tests were conducted on small-scale instances using two exact solution methods for multi-objective problems: weighted sum and augmented ε-constraint methods. The results obtained allowed us to explore the trade-off relationship between both objectives while evaluating the computational efficiency of both multi-objective solution methods.

**Resumen**

La compacidad de las rutas en los planes de distribución es un criterio que no ha sido lo suficientemente estudiado en la literatura relacionada a la distribución logística pero que ha demostrado tener un impacto relevante en la implementación práctica de los planes de ruteo, por ejemplo, en problemas de recolección de residuos sólidos urbanos. En este aspecto, este artículo presenta un modelo bi-objetivo para optimizar el problema de ruteo de vehículos con límites de tiempo considerando la minimización de los tiempos de viaje y la compacidad de las rutas. Se llevaron a cabo pruebas experimentales en instancias de pequeña escala utilizando dos métodos de resolución exacta para problemas multi-objetivo: la suma ponderada y el método de las ε-restricciones. Los resultados obtenidos permitieron explorar la relación de compromiso entre ambos objetivos a la vez que evaluar la eficiencia computacional de ambos métodos de resolución multiobjetivo.

**Resumo**

A compacidade das rotas em planos de distribuição é um critério que não foi suficientemente explorado na literatura relacionada à distribuição logística, mas tem demonstrado ter um impacto significativo na implementação prática de planos de roteamento, por exemplo, em problemas de coleta de resíduos sólidos urbanos. Nesse sentido, este artigo apresenta um modelo bi-objetivo para otimizar o problema de roteamento de veículos com restrições de tempo, considerando a minimização dos tempos de viagem e a capacidade das rotas. Testes experimentais foram realizados em instâncias de pequena escala usando dois métodos de resolução exata para problemas multi-objetivo: métodos de soma ponderada e ε-restrições aumentado. Os resultados obtidos nos permitiram explorar a relação de compromisso entre ambos os objetivos, ao mesmo tempo em que avaliamos a eficiência computacional de ambos os métodos de resolução multi-objetivo.

**Palabras clave:** Problemas de ruteo de vehículos; Compacidad; Optimización bi-objetivo.

**Palavras-chave:** Problemas de roteamento de veículos; Compacidade; Otimização bi-objetivo.

## 1. Introduction

In supply chain management, a significant part of the decision-making process involves logistics transportation planning. In the related literature, this type of planning is modeled through what is known as a Vehicle Routing Problem (VRP), a combinatorial optimization problem that aims at designing optimal routes for a fleet of vehicles to serve a set of customers (Toth and Vigo, 2014). When modeling a VRP, different objectives and constraints can be included depending on the specific application case (Arnold and Sörensen, 2019). In general, conventional or traditional optimization criteria in routing problems focus on objective measures such as minimizing transportation costs, distances or travel times, the number of required vehicles, or penalties for violating assumed "soft" constraints, such as time windows. Another optimization criteria that can be considered is the so-called "visual attractiveness" of routes, which is a more subjective measure. Visually appealing routes, primarily related to the degree of acceptance of solutions produced by VRP algorithms by route planners (Poot et al., 2002), have also demonstrated several advantages in the practical implementation and operational costs of the routing plans (Corberán et al., 2017). Although it is not easy to provide a precise definition of visual appealing routes due to the relatively subjectivity (Constantino et al., 2015), various authors associate this concept with some characteristics that the routes should have. These routes should be compact, confined to a bounded geographical region, and their convex hull has an approximately circular shape, should be separated from each other to avoid overlapping and they should not have intersections in their paths, preventing route crossovers (Rossit et al., 2019). This multidimensional nature of its definition has also led to a wide range of metrics developed to measure visual attractiveness in routing problems. An attempt to compare these various indicators and determine the degree of association between them was made in Rossit et al. (2019) through their application to well-known benchmarks in the literature.

Considering that the development of visually appealing distribution plans has been identified as a tool for reducing operational costs in transportation logistics, this work contributes with a bi-objective model to optimize both travel time and route compactness in a VRP with time limit. A computational experimentation is performed over realistic instances of the waste management system of Bahía Blanca to assess the trade-offs between the two objectives. The structure of this work is as follows. Section 2 presents the description of the problem, including the mathematical formulation. The experimentation and results are discussed in Section 3 and, finally, Section 4 presents the conclusions of the work and outlines future research directions.

## 2. The compact vehicle routing problem with time constraints

This section presents the compact vehicle routing problem with time constraints, including a conceptual description of the main features, related works and the mathematical formulation.

### 2.1. Conceptual problem and related works

Compactness is one of the most widely used measures to represent the visual attractiveness of a solution for a routing problem. However, despite being a relatively intuitive concept, compactness cannot be universally defined. Constantino et al. (2015) categorized the literature, distinguishing three types of compactness measures: i) similarity of the shape of the convex hull of the route with standard geometric shapes; ii) geographic/geometric or visual compactness; or iii) proximity between customers. In the same work, a different classification scheme was proposed, defining compactness measures based on: i) maximum travel times or Euclidean distances; ii) the sum of Euclidean distances; iii) the average and standard deviations of distances (or travel times) between customers and a reference point; or iv) the perimeters of zones or the perimeters and areas of zones. Among recent works applying compactness metrics, the study by Linfati et al. (2022) proposes various mathematical models to optimize two objectives in a weighted function: load balancing of routes and compactness. Kilby and Popescu (2021) introduced linear complexity algorithms for generating visually attractive routes. Tiwari and Sharma (2023) evaluated the performance of a set of commonly used heuristics from the literature, with one of the objectives in the aggregated cost function being compactness, and the tabu search procedure exhibited one of the best performances.

### 2.1 Mathematical formulation

Many of the compactness formulas developed in the literature (Rossit et al., 2019) can become complex when integrated into a mathematical formulation, leading to nonlinear models that are recognized as computationally challenging problems (Lee and Leyffer, 2011). Considering that vehicle routing problems are already NP-hard, even in their basic variants (Toth and Vigo, 2014), it is not surprising to see that most works incorporating compactness into optimization processes for routing problems are heuristic in nature. In this work, we propose a compactness calculation adapted from the work of Poot et al. (2002). So, the Vehicle Routing Problem (VRP) with a time limit that is being studied here, adapted for a municipal waste collection scenario, can be formalized through a mathematical programming model that includes the following definitions of sets, parameters, and variables.

| Sets | |
|---|---|
| $K$ | The set of vehicle routes that can be used. |
| $I$ | The set of points generating demands. |
| Parameters | |
| $C$ | The capacity of the vehicles. |
| $T$ | The time limit for the drivers' working day. |
| $c_i$ | The volume of the product to be collected at point $i$. |
| $t_{ij}$ | The travel time from point $i$ to point $j$ |
| $tc_i$ | The service time at point $i \in I$. |
| $td_0$ | The time for unloading the accumulated product at the depot. |

| $d_{ij}$ | The distance between nodes $i \in I$ and $j \in I$. |
|---|---|
| Variables | |
| $x_{ij}^k$ | A binary variable with a value of 1 if the vehicle travels from point $i \in I$ to point $j \in I$ during route $k \in K$, and 0 otherwise. |
| $z_{ij}^k$ | A binary variable with a value of 1 if nodes $i \in I$ and $j \in I$ are part of route $k \in K$. |
| $u_i^k$ | A continuous auxiliary variable introduced for the elimination of subtours. |
| $v_i^k$ | A continuous auxiliary variable introduced to limit the time per route. |

Taking into account these sets, parameters, and variables, the mathematical model for the problem can be expressed through the following Equations (1)-(10).

$$min \sum_{i,j \in M', \forall j \neq i, k \in K} t_{ij} x_{ij}^k \tag{1}$$

$$min \; COMP = \sum_{k \in K, i \in I, j \in I} d_{ij} z_{ij}^k \tag{2}$$

Subject to:

$$\sum_{i \in I, i \neq 0, j \neq i, k \in K} x_{ij}^k = 1, \qquad \forall j \in I \tag{3}$$

$$\sum_{i \in I, i \neq 0, j \neq i, k \in K} x_{ji}^k = 1, \qquad \forall j \in I \tag{4}$$

$$\sum_{j \in I, j \neq 0} x_{0j}^k \leq 1, \qquad \forall k \in K \tag{5}$$

$$\sum_{i \in I, j \neq i, k \in K} x_{ij}^k - \sum_{i \in I, j \neq i, k \in K} x_{ji}^k = 0, \qquad \forall j \in I \tag{6}$$

$$u_i^k - u_j^k \leq C(1 - x_{ij}^k) - c_j, \qquad \forall \; i, j \in I, j \neq i, k \in K \tag{7}$$

$$v_i^k - v_j^k \leq T(1 - x_{ij}^k) - t_{ij} - tc_j, \qquad \forall \; i, j \in I, j \neq i, j \neq 0, k \in K \tag{8}$$

$$v_i^k - v_j^k \leq T + T(1 - x_{ij}^k) - t_{i0} - td_0, \qquad \forall \; i, j \in I, j \neq i, k \in K \tag{9}$$

$$z_{ij}^k = \sum_{h \in I} x_{hi}^k * \sum_{h \in I} x_{hj}^k, \qquad \forall \; k \in K \tag{10}$$

$$\boldsymbol{x}, \boldsymbol{z} \in \mathcal{B} \quad \boldsymbol{u}, \boldsymbol{v} \geq 0$$

The objective functions are to minimize the total travel, defined by Equation (1), and minimize the compactness measure, defined by Equation (2). The model's constraints are presented in Equations (3) to (10). Equations (3) and (4) ensure that each node is visited only once, except for the depot (node 0). Equation (5) ensures that each route can start from the depot at most once. Equation (6) requires the conservation of flow within the route, while Equation (7) prevents the formation of subtours. Equations (8) and (9) ensure that the time per route does not exceed the maximum time established by the working day. Equation (10) allows $z_{ij}^k$ to be 1 if nodes $i$ and $j$ are part of the same route $k$ and 0 otherwise. The product of variables in Equation (10) is "non-linear". To linearize it, we apply the linearization methodology proposed by Glover (1984).

## 3. Computational experimentation

Preliminary computational experimentation was conducted on small sized test instances based on realistic instances of the municipal solid waste system of Bahía Blanca, Argentina that were developed in previous field work (Cavallin et al., 2020). The tests were performed using two different multi-objective approaches: the normalized weighted sum and the augmented ε-constraint method to study the different characteristics of both approaches (Rossit et al., 2017). The model was implemented using the Pyomo optimization package in a Python environment. The Gurobi solver version 10.0.2 was used for solving the resulting mathematical programming problem. The executions were carried out on a system with Intel® Xeon® Silver 4214 processors @ 2.2 GHz 2.19 GHz (2 processors), 32.0 GB of RAM, and a Windows® 10 Pro operating system.

In each instance, the number of efficient solutions to obtain is set to 10. For all the executions, a time limit of 1200 seconds was given to Gurobi to return the best solution obtained. In none of the executions Gurobi reached the optimal solution, which is not surprising given the NP-hard nature of VRP.

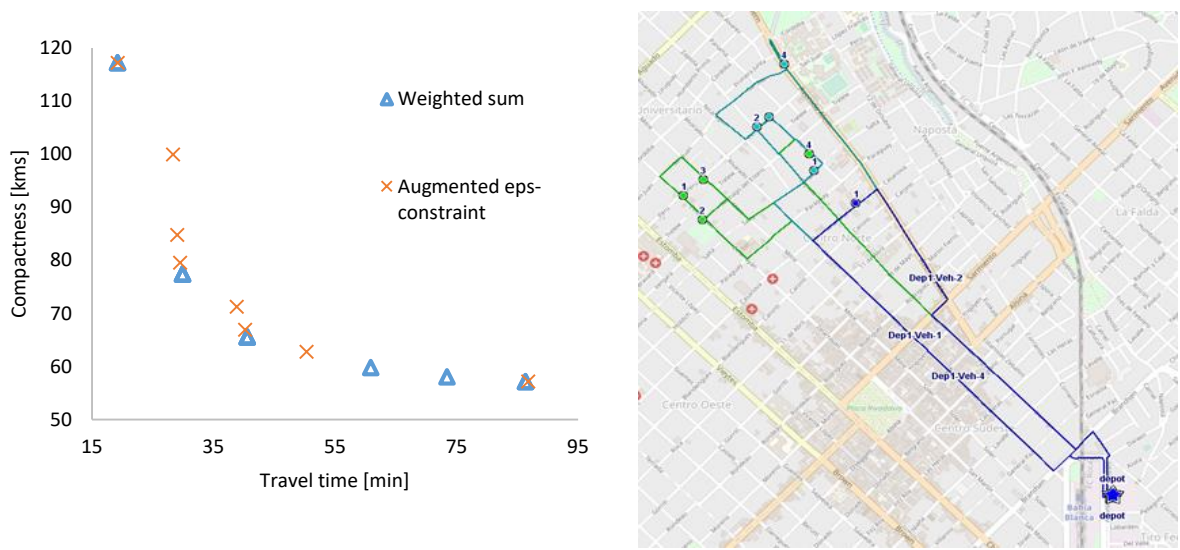

Figure 1. Pareto front of the proposed instance with 9 nodes (left) and non-dominated route of that front that provides a good trade-off between objectives (right).

For the instance with nine nodes, the augmented ε-constraint method was able to obtain seven non-dominated solutions in 8400 seconds since two of the runs were avoided due to the by-pass procedure of this method. However, given the simpler approach of the weighted sum, it does not have this capacity and performed the nine iterations and obtained only six non-dominated solutions in 12000 seconds since several solutions were repeated. A positive aspect of the weighted sum is that solutions are more evenly distributed along the Pareto front.

The results showed a clear trade-off among objectives, aiming at solutions that have a good performance in terms of travel time exhibit a bad results in terms of compactness and vice versa. This seems to confirm that using a bi-objective approach was a valid strategy for this problem. Additionally, the augmented ε-constraint method outperforms the weighting sum in terms of been able to find solutions of the non-convex regions of the Pareto front. Moreover, the augmented ε-constraint method was able to achieve a better usage of the computational resources.

## 4. Discussion and future work

A constant demand in current socioeconomic systems is finding new strategies for generating efficient solutions to problems related to transportation and distribution logistics while incorporating functionally appropriate performance criteria. This work focusesas on solving vehicle routing problems considering capacity and total route time constraints, with the objective of minimizing travel time and optimizing compactness of the obtained routes. Although only a preliminary computational experimentation over test instances is presented, the main contribution of this work is to present a bi-objective mathematical formulation for the VRP with time limits that can be linearized. As aforementioned, these kinds of problems have not been sufficiently studied in the related literature. The results for the test instances have shown a clear trade-off between both objectives, highlighting that minimizing travel times and route compactness are conflicting goals when evaluating the generated solutions. In this sense, providing the decision-maker with a set of non-dominated solutions is beneficial because it allows the selection of the most suitable solution based on the specific scenario, seeking a reasonable compromise between travel times and route compactness. Moreover, this work compared two traditional exact multiobjective approaches: the weighted sum and the ε-constraint methods. The ε-constraint method surpasses the weighting sum in terms of been able to find solutions of non-convex regions of the Pareto front. Also, the augmented ε-constraint method was able to achieve a better usage of the computational resources.

Future research lines can be oriented toward expanding the computational experimentation with more diverse instances in which the objective of compactness can be better assessed. Additionally, heuristic approaches can be developed to addressed the trade-off between travel times and compactness in real-world (larger) instances.

### References

ARNOLD, F., & SÖRENSEN, K. (2019). What makes a VRP solution good? The generation of problem-specific knowledge for heuristics. *Computers & Operations Research*, 106, 280-288.

CAVALLIN, A., ROSSIT, D., HERRAN, V., ROSSIT, D., & FRUTOS, M. (2020). Application of a methodology to design a municipal waste pre-collection network in real scenarios. *Waste Management & Research*, 38(1_suppl), 117-129.

CONSTANTINO, M., GOUVEIA, L., MOURÃO, M., & NUNES, A. (2015). The mixed capacitated arc routing problem with non-overlapping routes. *European Journal of Operational Research*, 244(2), 445-456.

CORBERÁN, Á., GOLDEN, B., LUM, O., PLANA, I., & SANCHIS, J. (2017). Aesthetic considerations for the min-max K-Windy Rural Postman Problem. *Networks*, 70(3), 216-232.

KILBY, P., & POPESCU, D. (2021). Linear Complexity Algorithms for Visually Appealing Routes in the Vehicle Routing Problem. In *Data and Decision Sciences in Action 2: Proceedings of the ASOR/DORS Conference 2018* (pp. 81-98). Springer International Publishing.

LEE, J., & LEYFFER, S. (Eds.). (2011). Mixed integer nonlinear programming (Vol. 154). Springer Science & Business Media.

LINFATI, R., YÁÑEZ, F., & ESCOBAR, J. (2022). Mathematical models for the vehicle routing problem by considering balancing load and customer compactness. *Sustainability*, 14(19), 12937.

POOT, A., KANT, G., & WAGELMANS, A. P. M. (2002). A savings based method for real-life vehicle routing problems. *Journal of the Operational Research Society*, 53(1), 57-68.

ROSSIT, D., TOHMÉ, F., FRUTOS, M., & BROZ, D. (2017). An application of the augmented ε-constraint method to design a municipal sorted waste collection system. *Decision Science Letters*, 6(4), 323-336.

ROSSIT, D., VIGO, D., TOHMÉ, F., & FRUTOS, M. (2019). Visual attractiveness in routing problems: A review. *Computers & Operations Research*, 103, 13-34.

TIWARI, K., & SHARMA, S. (2023). An optimization model for vehicle routing problem in last-mile delivery. *Expert Systems with Applications*, 222, 119789.

TOTH, P. & VIGO, D. (2014) The vehicle routing problem. Society for Industrial and Applied Mathematics, Philadelphia, United States of America.